\documentclass{article}
\usepackage{amsmath,amssymb,amsfonts,amsthm}
\makeindex
\title{\textbf{Independently Axiomatizable $\lomegaone$ Theories}    \\ }
\author{\begin{tabular}{c}Greg Hjorth\\Department of Mathematics,
\\The University of Melbourne, \\ Parkville, VIC 3010, Australia\\
 G.Hjorth@ms.unimelb.edu.au\\\end{tabular}
\begin{tabular}{c}Ioannis A. Souldatos\\
273 Wissink Hall, Office 261\\
Mathematics Department\\
Minnesota State University\\
Mankato, MN 56001\\ioannis.souldatos@mnsu.edu\\\end{tabular}
  }
\begin{document}
\newcommand{\omegaone}{\ensuremath{\omega_1}}
\newcommand{\lomegaone}{\ensuremath{\mathcal{L}_{\omega_1,\omega}}}
\newcommand{\alephs}[1]{\ensuremath{\aleph_{#1}}}
\newcommand{\alephalpha}{\alephs{\alpha}}
\newcommand{\alephomegaone}{\alephs{\omegaone}}
\newcommand{\alephalphaplus}{\alephs{\alpha+1}}
\newcommand{\atypes}{\ensuremath{\alpha\mbox{-types}}}
\newcommand{\atype}{\ensuremath{\alpha\mbox{-type}}}
\newcommand{\btypes}{\ensuremath{\beta\mbox{-types}}}
\newcommand{\btype}{\ensuremath{\beta\mbox{-type}}}
\newcommand{\z}{\ensuremath{\mathcal{Z}}}
\newcommand{\M}{\ensuremath{\mathcal{M}}}
\newcommand{\N}{\ensuremath{\mathcal{N}}}
\newcommand{\A}{\ensuremath{\mathcal{A}}}
\newcommand{\B}{\ensuremath{\mathcal{B}}}
\newcommand{\F}{\ensuremath{\mathcal{F}}}
\newcommand{\D}{\ensuremath{\mathcal{D}}}
\newcommand{\C}{\ensuremath{\mathcal{C}}}
\newcommand{\E}{\ensuremath{\mathcal{E}}}
\newcommand{\G}{\ensuremath{\mathcal{G}}}
\newcommand{\V}{\ensuremath{\mathbb{V}}}
\newcommand{\Vprime}{\ensuremath{\mathbb{V'}}}
\newcommand{\veca}{\ensuremath{\vec{a}}}
\newcommand{\vecb}{\ensuremath{\vec{b}}}
\newcommand{\vecc}{\ensuremath{\vec{c}}}
\newcommand{\vecd}{\ensuremath{\vec{d}}}
\newcommand{\vecx}{\ensuremath{\vec{x}}}
\newcommand{\vecy}{\ensuremath{\vec{y}}}
\newcommand{\vecz}{\ensuremath{\vec{z}}}
\newcommand{\lang}[1]{\ensuremath{\mathcal{L}_{#1}}}
\newcommand{\langhat}{\ensuremath{\widehat{\lang{}}}}
\newcommand{\phialpha}{\ensuremath{\phi_\alpha}}
\newcommand{\phibeta}{\ensuremath{\phi_\beta}}
\newcommand{\psialpha}{\ensuremath{\psi_\alpha}}
\newcommand{\psibeta}{\ensuremath{\psi_\beta}}
\newcommand{\phikappa}{\ensuremath{\phi_\kappa}}
\newcommand{\phiM}{\ensuremath{\phi_{\M}}}
\newcommand{\philtok}{\ensuremath{\phi_{\ltok}}}
\newcommand{\phialphaplus}{\ensuremath{\phi_{\alpha+1}}}
\newcommand{\phibaralpha}{\ensuremath{\overline{\phi_\alpha}}}
\newcommand{\phibarbeta}{\ensuremath{\overline{\phi_\beta}}}
\newcommand{\phialphazero}{\ensuremath{\phi^{(\alpha)}_0}}
\newcommand{\phialphan}{\ensuremath{\phi^{(\alpha)}_n}}
\newcommand{\phiplus}{\ensuremath{\phi^+}}
\newcommand{\phistar}{\ensuremath{\phi^*}}
\newcommand{\phivecaMalpha}{\ensuremath{\phi^{vec{a},\M}_{\alpha}}}
\newcommand{\phipar}[2]{\ensuremath{\phi^{#1}_{#2}}}
\newcommand{\phialphapar}[1]{\ensuremath{\phi^{#1}_{\alpha}}}
\newcommand{\phibetapar}[1]{\ensuremath{\phi^{#1}_{\beta}}}
\newcommand{\scM}{\ensuremath{\alpha(\M)}}
\newcommand{\Malpha}{\ensuremath{\mathcal{M}_\alpha}}
\newcommand{\continuum}{\ensuremath{2^{\alephs{0}}}}
\newcommand{\Talpha}{T_{\alpha}}
\newcommand{\Tbeta}{T_{\beta}}
\newcommand{\Tstar}{T^{\star}}
\newcommand{\Sinfty}{\ensuremath{S_{\infty}}}
\newcommand{\f}[1]{\ensuremath{f(\{#1\})}}
\newtheorem{thrm}{Theorem}
\newtheorem{lem}[thrm]{Lemma}
\newtheorem{cor}[thrm]{Corollary}
\newtheorem{df}[thrm]{Definition}
\newtheorem{claim}{Claim}
\newtheorem{note}[thrm]{Note}
\maketitle
Note: This paper has been published by the Journal of Symbolic Logic which owns the copyright. 

Bibliographical Reference:  J. Symbolic Logic, Volume 74, Issue 4 (2009), 1273-1286 

Permanent Link: http://projecteuclid.org/euclid.jsl/1254748691

\abstract In partial answer to a question posed by Arnie Miller \cite{Miller'sWebpage} and X. Caicedo \cite{Caicedo} we obtain sufficient conditions for an $\lomegaone$ theory to have an independent axiomatization. As a consequence we obtain two corollaries: The first, assuming Vaught's Conjecture, every $\lomegaone$ theory in a countable language has an independent axiomatization. The second, this time outright in ZFC, every intersection of a family of Borel sets can be formed as the intersection of a family of \emph{independent} Borel sets.

\section{Introduction}\begin{df}\label{IndDef} A set of sentences $T'$ is called \emph{independent} if for every $\phi\in T'$, $T'\setminus \{\phi\}\nvDash\phi$. A theory $T$ is called \emph{independently axiomatizable}, if there is a set $T'$ which is independent and $T$ and $T'$ have exactly the same models.\end{df}

Note that this definition applies to sets of sentences in both first-order ($\lang{\omega,\omega}$) and infinitary ($\lomegaone$) logic, granted that we have defined a meaning for $\models$. The question is whether every theory has an independent axiomatization. For first-order theories the answer is positive:

\begin{thrm} (M.I. Reznikoff- \cite{Reznikoff}) All theories of any cardinality in $\lang{\omega,\omega}$, are independently axiomatizable.
\end{thrm}

\begin{df}\label{ImplyDef} For a set of sentences $T\subset\lomegaone$ and a sentence $\sigma\in\lomegaone$, write \[T\models\sigma,\] if all the models of $T$ satisfy $\sigma$. Two sets of sentences in $\lomegaone$ are \emph{semantically equivalent} if they have exactly the same models.\end{df}

Using definition \ref{IndDef} and \ref{ImplyDef}, it makes sense to ask whether a theory in $\lomegaone$ is independently axiomatizable, i.e. when it is semantically equivalent to an independent set. A partial result to this question is given by
\begin{thrm}\label{CaicedosThrm} (X. Caicedo- \cite{Caicedo}) Any theory in$\lomegaone$ of cardinality no more than $\aleph_1$ is independently axiomatizable. \end{thrm}For cardinalities greater than $\aleph_1$, Caicedo obtained partial results for a weaker notion of \emph{countable independence}, which requires that every countable subset of the set of sentences is independent. Our main result (theorem \ref{MainTheorem})  states that\begin{thrm}\label{MainThrmIntro} For a countable language $\lang{}$ and for a theory $T\subset\lomegaone$, if the number of counterexamples to Vaught's Conjecture contained in $T$ is \emph{small}, then $T$ is independently axiomatizable.\end{thrm}

The meaning of a \emph{small} number of counterexamples is made clear by definition \ref{SmallDef}. Vaught's conjecturestates\vspace{0.2cm}

\textbf{Conjecture} (Vaught) Every sentence $\sigma\in\lomegaone$ either has countable many non-isomorphic countable models, or else it has continuum many. \vspace{0.2cm}

Under the Continuum Hypothesis the conjecture is trivially true and Morley proved that every counterexample to it will necessarily have $\aleph_1$ many non-isomorphic countable models. Using theorem \ref{MainThrmIntro} we then obtain two consequences:

\begin{thrm}Assume Vaught's Conjecture. Then for any countable $\lang{}$ and any theory $T\subset\lomegaone$, $T$ is independently axiomatizable.\end{thrm}

\begin{df}We call a collection of Borel sets $\B=\{B_i|i\in I\}$ \emph{independent} if $\bigcap\B\neq\emptyset$ and for every $i\in I$, \[\bigcap_{j\neq i} B_j\setminus B_i\neq\emptyset.\]Two collections $\B,\B'$ are \emph{equivalent}, if\[\bigcap \B=\bigcap\B'.\]\end{df}

\begin{thrm} Let ${\mathcal B}$ be a collection of Borel sets. Then there is an independent collection of Borel sets ${\mathcal B'}$ with\[\bigcap {\mathcal B}=\bigcap{\mathcal B'}.\]\end{thrm}

\section{Preliminary work}

\hspace{.45cm}Our results are about sentences in \lomegaone. When we refer to a sentence, we mean a sentence in $\lomegaone$ and when we refer to a theory, we mean a theory of sentences in \lomegaone. If a theory $T$ doesn't have any models, it is axiomatizable by the sentence $\exists x (x\neq x)$, while any theory $T$ is equivalent to the theory obtained by deleting from $T$ all its valid sentences. Thus, we can assume that all the theories we work with are consistent and do not contain valid formulas. Throughout this paper  we assume that the language $\lang{}$ we are working with is countable. Then every theory $T\subset\lomegaone$ can have size up to the continuum. Under the Continuum Hypothesis and using theorem \ref{CaicedosThrm}, every theory is independently axiomatizable. So, it suffices to deal with the case that the Continuum Hypothesis fails and we will take this as one of our working assumptions.

\begin{df} Let $\phi$ a sentence. We say that the sentences $\{\psialpha|\alpha\in I \}$ \emph{partition} $\phi$ if:
\begin{itemize}
\item for all $\alpha$, $\psialpha$ is consistent,
\item $\models \phi \leftrightarrow \bigvee_{\alpha\in I} \psialpha$ and
\item for all $\alpha$, $\models \psialpha \rightarrow     \bigwedge_{\beta\neq\alpha} \neg\psi_{\beta}$
\end{itemize}
\end{df}
What we are heading towards is to prove that under the assumption of a perfect set of countable models, we get a partition into continuum many sentences.
If $\M$ is a countable model and $\vec{a}\in\M$, define the \emph{$\alpha$-type} of $\vec{a}$ in $\M$ inductively:

\begin{eqnarray*}  \phi^{\vec{a},\M}_0&:=& \bigwedge\{\psi(\vec{x})|\psi\mbox{ is an atomic formula or negation of atomic},\M\models\psi(\vec{a})\}, \\  \phi^{\vec{a},\M}_{\alpha+1} &:=&  \phi^{\vec{a},\M}_{\alpha}\bigwedge\{\exists\vec{y}\phi^{\vec{a}\frown\vec{b},\M}_{\alpha}(\vec{x},\vec{y})|\vec{b}\in\M\}\wedge \\  & & \bigwedge_n\forall y_0\ldots y_n\bigvee\{\phi^{\vec{a}\frown\vec{b},\M}_{\alpha}(\vec{x},\vec{y})|\vec{b}\in\M\},\\  \phi^{\vec{a},\M}_{\lambda} &:=&  \bigwedge_{\alpha<\lambda}\phi^{\vec{a},\M}_{\alpha}\mbox{, for $\lambda$ limit.}
\end{eqnarray*}

The \emph{$\alpha$-types of $\M$} are defined to be all sentences of the form $\phi^{\vec{a},\M}_{\alpha}$, for any $\vec{a}\in\M$, and if $\sigma$ is a sentence, the \emph{$\alpha$-types of $\sigma$} are all sentences of the form $\phi^{\vec{a},\M}_{\alpha}$ with $\M\models\sigma$ and $\vec{a}\in\M$. If $\M$ is a countable model, then it realizes only countably many types and there is an ordinal $\delta<\omegaone$ such that for all $\vec{a},\vec{b}\in\M$,

\[\phi^{\vec{a},\M}_{\delta}\neq\phi^{\vec{b},\M}_{\delta}\mbox{ iff there is some $\gamma<\omegaone$ }(\phi^{\vec{a},\M}_{\gamma}\neq \phi^{\vec{b},\M}_{\gamma}).\]

The least such ordinal $\delta$ we call the Scott height of $\M$ and write $\scM$. Then $\phipar{\emptyset,\M}{\scM+2}$ is called the Scott sentence of $\M$.

\begin{df} For a $\lomegaone$-sentence $\phi$ and $\alpha<\omegaone$,  let\[\Psi_{\alpha}(\phi):=\{\phi^{\vec{a},\M}_\alpha|\vec{a}\in\M,\M\models \phi\},\]the $\atypes$ of $\phi$. Define also\[\Phi_{\alpha}(\phi):=\{\phi^{\emptyset,\M}_\alpha|\M\models \phi\}.\]
\end{df}

Now, observe that if $\alpha=\gamma+1$, some $\gamma$, then we can identity every $\phi^{\vec{a},\M}_{\gamma+1}$ with the set\[\{\phi^{\vec{a}\frown\vec{b},\M}_{\gamma}|\vec{b}\in\M\}.\]

This enables us to consider $\Psi_{\gamma+1}(\phi)$ and $\Phi_{\gamma+1}(\phi)$ as  subsets of $X_\alpha(\phi):=2^{\Psi_{\gamma}(\phi)}$. In the special case that $\Psi_{\gamma}(\phi)$ is countable, $X_\alpha(\phi)$ becomes a standard Borel space and we will prove (lemma\ref{someborelsets}) that in this case $\Psi_{\gamma+1}(\phi)$ and $\Phi_{\gamma+1}(\phi)$ are $\mathbf{\Sigma^1_1}$ subsets.
Similarly, for $\alpha$ limit, we can identify $\phi^{\vec{a},\M}_{\alpha}$ with the set\[\{\phi^{\vec{a},\M}_{\gamma}|\gamma<\alpha\}.\]
Then $\Psi_{\alpha}(\phi)$ and $\Phi_{\alpha}(\phi)$ become subsets of $X_\alpha(\phi):=2^{\bigcup_{\gamma<\alpha} \Psi_{\gamma}(\phi)}$. Again, in the case that for all $\gamma<\alpha$, $\Psi_{\gamma}(\phi)$ is countable, $X_\alpha(\phi)$ becomes a standard Borel space and $\Psi_{\alpha}(\phi)$ and $\Phi_{\alpha}(\phi)$ are $\mathbf{\Sigma^1_1}$ subsets. The same can be said for $\phi^{\vec{a},\M}_0$. We can identify it with \[\{\phi(x)|\mbox{$\phi$ atomic},\M\models\phi(\vec{a})\},\] in which case $\Psi_0(\phi)$ and $\Phi_0(\phi)$ become subsets of $X_0(\phi):=2^A$, with $A$ being the set of all atomic, or negation of atomic sentences. Since we assumed that the language we work with is countable, $A$ is countable and $X_0(\phi)$ is a standard Borel space with $\Psi_0(\phi)$ and $\Phi_0(\phi)$ $\mathbf{\Sigma^1_1}$ subsets.

\begin{df} Let $\lang{}$ be a countable language and let $Mod(\lang{})$ be the set of all countable $\lang{}$-structures with underlying set $\mathbb{N}$. We equip $Mod(\lang{})$ with the topology generated by taking as basic open sets all sets of the form
\[\{M\in Mod(\lang{})| M\models \varphi(n_1,...n_m)\}\]
for $\varphi(\vec x)$ a quantifier free formula and $n_1,...,n_m \in\mathbb{N}$.
\end{df}
It is easily shown that $Mod(\lang{})$ is a Polish space. For more on this one can consult \cite{BeckerKechrisBook}.

\begin{df} For a sentence $\sigma$ let $Mod(\sigma)$ be the set of all models in $Mod(\lang{})$ that satisfy $\sigma$.
\end{df}

This becomes a standard Borel space space by the Borel structure it inherits from $Mod(\lang{})$ (cf. \cite{BeckerKechrisBook} too).

\begin{lem}\label{someborelsets} Let $\phi$ be a $\lomegaone$-sentence, $\alpha<\omegaone$, $\Psi_{\alpha}(\phi)$, $\Phi_{\alpha}(\phi)$ and $X_\alpha(\phi)$ as defined above. Assume that for all $\gamma<\alpha$, $\Psi_{\gamma}(\phi)$ is countable. Then
\begin{enumerate}
\item the function $Mod(\phi)\times \omega^{<\omega}\rightarrow  X_\alpha(\phi)$,
with \[(\M,\vec{a})\mapsto \phi^{\vec{a},\M}_{\alpha}\]  is Borel and
\item $\Psi_{\alpha}(\phi)$ and $\Phi_{\alpha}(\phi)$ are   $\mathbf{\Sigma^1_1}$ sets.
\end{enumerate}

\begin{proof} Recall that under the countability assumption for the $\Psi_{\gamma}$'s, $X_\alpha(\phi)$ becomes a standard Borel space with $\Psi_{\alpha}(\phi)$ and $\Phi_{\alpha}(\phi)$ seen as subsets of it. Therefore, the statement of the theorem makes sense. Now, by induction on $\beta\le\alpha$, it follows easily from the definition that the function $(\M,\vec{a})\mapsto\phi^{\vec{a},\M}_\beta$ is Borel. In particular, the same is true for $(\M,\vec{a})\mapsto \phi^{\vec{a},\M}_\alpha$. Using this function we can write
\[\psi\in\Psi_{\alpha}\mbox{ iff }\exists\M\exists\vec{a}\in\M((\M\models\phi) \mbox{ and } (\psi=\phi^{\vec{a},\M}_\alpha)),\]
and similarly
\[\psi\in\Phi_{\alpha}\mbox{ iff }\exists\M((\M\models\phi) \mbox{ and } (\psi=\phi^{\emptyset,\M}_\alpha)).\]

This proves the lemma.\end{proof}
\end{lem}
If $\Psi_\alpha(\phi)$ is as in the above lemma, then by the perfect set theorem for $\mathbf{\Sigma^1_1}$ sets, it is either countable or has size continuum. If it is countable, then we can apply the lemma once more and we can keep doing that until we either run out of countable ordinals, or until we find an uncountable $\Psi_{\alpha'}(\phi)$, some $\alpha'>\alpha$.

\begin{lem}\label{ScottHeight} For all $\alpha<\omegaone$, the set
\[\{(\M,\N,\veca,\vecb)\in Mod(\lang{})^2\times (\omega^{<\omega})^2|\phialphapar{\veca,\M}=\phialphapar{\vecb,\N}\}\]
is Borel. In particular, for $\phi\in\lomegaone$ and $\gamma<\omegaone$, the set\[\{\M\in Mod(\phi)|\alpha(\M)<\gamma\}\] is also Borel.\begin{proof} For the first part, by induction on $\alpha$:

$\alpha=0:$ Then $\phipar{\veca,\M}{0}=\phipar{\vecb,\N}{0}$ if and only if for every atomic, or negation of atomic, formula $\phi$,
\[\M\models\phi(\veca)\leftrightarrow\N\models\phi(\vecb).\]

$\alpha+1:$ Then $\phipar{\veca,\M}{\alpha+1}=\phipar{\vecb,\N}{\alpha+1}$ if and only if \[\forall \vecc\in\M\exists\vecd\in\N(\phipar{\veca\frown\vecc,\M}{\alpha}=
\phipar{\vecb\frown\vecd,\N}{\alpha})\]
and
\[\forall\vecd\in\N \exists\vecc\in\M(\phipar{\veca\frown\vecc,\M}{\alpha} =\phipar{\vecb\frown\vecd,\N}{\alpha}).\]

$\alpha$ limit: Then $\phipar{\veca,\M}{\alpha}=\phipar{\vecb,\N}{\alpha}$ if and only if \[\forall \beta<\alpha(\phipar{\veca,\M}{\beta}=\phipar{\vecb,\N}{\beta}).\]

By inductive hypothesis, all these conditions are Borel and therefore our set is Borel. Now, by the definition of the Scott height, $\scM<\gamma$ if and only if \[\bigvee_{\alpha<\gamma}\forall\veca\;\forall\vecb\;(\phipar{\veca,\M}{\alpha}
=\phipar{\vecb,\M}{\alpha}\rightarrow\phipar{\veca,\M}{\alpha+1}
=\phipar{\vecb,\M}{\alpha+1}).\]

By the first part, this condition is Borel.
\end{proof}
\end{lem}

\begin{lem} If a $\lomegaone$-sentence $\phi$ has continuum many non-isomorphic countable models, then there are countable ordinals $\alpha<\beta$, with $X_{\alpha}(\phi)$ a standard Borel space, a perfect set $P$ and continuous functions $t:P\rightarrow X_{\alpha}(\phi)$, $M:P\rightarrow Mod(\phi)$ such that:
    \begin{itemize}
     \item for all $x\neq y\in P$, $t(x),t(y)$ are distinct types in $X_{\alpha}(\phi)$,
     \item for all $x\in P$, $M(x)$ is a countable model of $\phi$ that realizes    $t(x)$ and has Scott height $<\beta$. 
     \end{itemize}

     Moreover, we can assume that for $x\neq y\in P$, $M(x)\nvDash t(y)$.
     \begin{proof} Let $\alpha<\omegaone$ be the least ordinal with $\Phi_{\alpha}(\phi)$ uncountable. Then $\Phi_{\gamma}(\phi)$, $\gamma<\alpha$, are all countable and applying lemma \ref{someborelsets}, we conclude that $X_{\alpha}(\phi)$ is a standard Borel space and $\Phi_{\alpha}(\phi)$ is $\mathbf{\Sigma^1_1}$. Consider an ordinal $\beta>\alpha$ large enough so that the set\[\{\phi^{\emptyset,\M}_{\alpha}\in\Phi_{\alpha}(\phi)|\alpha(\M)<\beta\}\]is still uncountable. By lemma \ref{ScottHeight} this is again $\mathbf{\Sigma^1_1}$. Consequently, it embeds a perfect set. So, let $P$ a perfect set with \[t:P\rightarrow\{\phi^{\emptyset,\M}_{\alpha}\in\Phi_{\alpha}(\phi)|\alpha(\M)<\beta\}\subset X_{\alpha}(\phi)\] a continuous embedding. Then every $t(x)$, $x\in P$, has the form $\phi^{\emptyset,\M}_\alpha$, for some $\M$ with $\alpha(\M)<\beta$. Consider the set
     \[\{(x,\M)\in P\times Mod(\phi)|\M\models t(x),\alpha(\M)<\beta\}.\]%
     This is not empty and by lemma \ref{ScottHeight} and since $t$ is continuous, it is Borel. By Jankov- von Neumann Uniformization theorem (cf. \cite{Kechrisbook}), we get a function $x\mapsto M(x)$ that is Baire measurable and for all $x\in P$, $(x,M(x))$ is in the above set. Restricting the domain to a comeager set $C\subset P$ we can further assume that $x\mapsto M(x)$ is continuous on $C$. Let \[R_0:=\{(x,y)\in C^2|M(x)\nvDash t(y)\}.\] Since $t$ is 1-1 and $M(x)$ can satisfy only countably many $\atypes$, $R_0$ is comeager in $C^2$. By theorem 19.1 of \cite{Kechrisbook}, we can find a Cantor set $C_1\subset C$ such that \[(C_1)^2\subset R_0,\] and then for all $x,y\in C_1$, \[x\neq y\Rightarrow (x,y)\in R_0\Rightarrow M(x)\nvDash t(y),\] which proves the lemma.

\end{proof}
\end{lem}

Note: The previous lemma is the only place where we use the assumption $\alephs{1}<2^{\alephs{0}}$.

Observe also that for $x\neq y\in P$, $M(x)\models t(x)$, while $M(y)\nvDash t(x)$, which implies that $M(x)\ncong M(y)$.

\begin{lem} The set $A_0:=\{\M|\exists x\in P (\M\cong M(x))\}$ is Borel.
\begin{proof} We need first that the set $A_1:=\{(x,\M)|\M\cong M(x),x\in P\}$ is Borel. Since for all $x\in P$ the Scott height of $M(x)$ is $<\beta$,
\begin{eqnarray*}  \M\cong M(x) &iff& \M\models\phipar{\emptyset,M(x)}{\beta+1} \\
   &iff&   \phipar{\emptyset,\M}{\beta+1}=\phipar{\emptyset,M(x)}{\beta+1}.
\end{eqnarray*}
This last condition is Borel by lemma \ref{ScottHeight}. By the observation that for $x\neq y\in P$, $M(x)\ncong M(y)$, we can also conclude that if $(x_1,\M), (x_2,\M)$ are both in $P$, then $x_1=x_2$. By the Lusin-Novikov theorem, the projection of $A_1$ (on the second component) is also Borel and this is exactly what we have to prove.
\end{proof}
\end{lem}

\begin{cor} There is a sentence $\phiplus\in\lomegaone$ such that for every countable model $\M$, \[\M\models\phiplus\mbox{ iff }\M\in A_0.\]
\begin{proof} $A_0$ is obviously invariant under isomorphisms and by the previous lemma is Borel. Therefore,there exists a $\lomegaone$ sentence $\phiplus$ as in the statement.
\end{proof}
\end{cor}

\begin{lem} If $\N$ is a model of $\phiplus$, countable or uncountable, and it satisfies one of the $\atypes$ $\{t(x)|x\in P\}$, then it actually satisfies the Scott sentence $s(x)$ of $M(x)$.
\begin{proof} Recall here that a countable model can satisfy only one of the $\atypes$ $\{t(x)|x\in P\}$. If $\N$ is countable and satisfies $\phiplus$, then it belongs to$A_0$, i.e. it is isomorphic to one of the $M(x)$, $x\in P$. If $\N\models t(x)$, then $M(x)\cong \N$ and obviously $\N\models s(x)$. Therefore, assume that $\N$ is uncountable with $\N\models t(x)$, some $x\in P$. Let $s(x)$ be the Scott sentence of $M(x)$ and $\F$ the fragment generated by $\phiplus,t(x)$ and $s(x)$. Let $\N_0$ be a countable model with \[\N_0\prec_{\F}\N.\] Then $\N_0\models \phiplus$ and $\N_0\models t(x)$. As in the countable case, $\N_0\models s(x)$, which implies that $\N\models s(x)$.
\end{proof}
\end{lem}

Using all these lemmas we are ready to prove

\begin{thrm}\label{partitionphi} If $\phi$ has $2^{\aleph_0}$ many non-isomorphic countable models, then there exists a partition of $\phi$ into continuum many sentences.
\begin{proof} Assume that $P,\alpha,x\mapsto t(x),x\mapsto M(x)$ and $\phiplus$ are as above.
\begin{claim} It suffices to find a $\lomegaone$-sentence $\phistar$ that expresses the fact that our model satisfies one of the $\atypes$ $\{t(x)|x\in P\}$.
\begin{proof}(of claim) First we need that every model of $\phiplus$ is also a model of $\phistar$. Arguing as before let $\N\models\phiplus$, $\F$ be the fragment generated by both $\phiplus$ and $\phistar$, and $\N_0\prec_{\F}\N$ a countable model. Then, there exists $x\in P$ with $\N_0\cong M(x)$ and $\N_0\models t(x)$. By definition, $\N_0\models\phistar$ which also implies that $\N\models\phistar$. Combining this with the previous lemma, we conclude that every countable or uncountable model of $\phiplus$ will satisfy one of the Scott sentences $\{s(x)|x\in P\}$. Therefore,
\[\{\phi\wedge\neg\phiplus\}\cup\{s(x)|x\in P,\mbox{ $s(x)$ is the Scott sentence of $M(x)$}\}\]
gives a partition of $\phi$ into continuum many sentences.
\end{proof}
\end{claim}

Towards constructing $\phistar$, let $S:=\bigcup_{\gamma<\alpha}\Psi_{\gamma}(\phi)$. By assumption on $\alpha$, $S$ is countable and for all $x\in P$, $t(x)\in X_\alpha(\phi)\subset 2^S$. For all $u\in 2^{<\omega}$ we can construct $S_u$ finite subsets of $S$ such that\begin{enumerate}  \item for every $u$, $S_{u\frown 0}$ is always incompatible with  $S_{u\frown 1}$,  \item $S_u\subset S_w$ when $u\subset w$,  \item for every $\hat{u}\in 2^{\omega}$,  $\bigcup_{n\in \mathbb{N}} S_{\hat{u}\upharpoonright n}$ is an  element of $\{t(x)|x\in P\}$ and every $t(x)$ in this set can be  written as $\bigcup_{n\in \mathbb{N}} S_{\hat{u}\upharpoonright n}$,  for some $\hat{u}\in 2^{\omega}$.\end{enumerate}Consider the sentence:\[\phistar:=\exists a\bigwedge_{n\in\omega} \bigvee_{u\in 2^n}\bigwedge_{\psi\in S_u} \psi(a).\]It is obvious that every model of $\phistar$ will satisfy one of the $\atypes$ $t(x)$, $x\in P$.
\end{proof}
\end{thrm}

\section{Main Result}

\hspace{.5cm}We work as before with a countable $\lang{}$. Throughout this section we will not distinguish between a model $\M$ and its isomorphism class $[\M]_{\Sinfty}$. So, when we say that a sentence has countably many countable models, we actually mean countably many non-isomorphic countable models.

\begin{df}\label{t0t1t2def} For a theory $T=\{\phi_\alpha|\alpha<2^{\aleph_0}\}$ define
\begin{eqnarray*}
T_0 &:=& \{\phi\in T|\neg\phi\mbox{ has countably many countable models}\}, \\
T_1 &:=& \{\phi\in T|\neg\phi\mbox{ has $\aleph_1$many countable models}\}, \\
T_2 &:=& \{\phi\in T|\neg\phi\mbox{ has $2^{\aleph_0}$ many countable models}\},
\end{eqnarray*}
and
\begin{eqnarray*}
X_0(T) &:=& \{\M|\M\models\neg\phi\mbox{, some } \phi\in T_0, \mbox{ $\M$ countable}\} \\  X_1(T) &:=& \{\M|\M\models\neg\phi\mbox{, some } \phi\in T_1, \mbox{ $\M$ countable}\} \\  X_2(T) &:=& \{\M|\M\models\neg\phi\mbox{, some } \phi\in T_2, \mbox{ $\M$ countable}\}\\  X(T) &:=& X_0(T)\cup X_1(T)\cup X_2(T).
\end{eqnarray*}
\end{df}

Note that the sets $T_0, T_1$ and $T_2$ are disjoint, while the sets $X_0(T), X_1(T)$ and $X_2(T)$ may not be disjoint. Also, all sentences in $T_1$ provide counterexamples to Vaught's Conjecture.

\begin{df}\label{SmallDef} In case that $|X(T)|\ge|T_1|$ we will say that $T_1$ is \emph{small} in $T$, or if $T$ is apparent, we will just say that $T_1$ is \emph{small}.
\end{df}

Smallness assumption for $T_1$ will be crucial for our main result (theorem \ref{MainTheorem}). If $|X(T)|=\aleph_0$, then $T_1=T_2=\emptyset$ and if $\{\M_n|n\in\omega\}$ enumerate the models in $X_0$ and $\{\phi_n|n\in\omega\}$ enumerate their Scott sentences, then it is easily seen that $T$ is equivalent to
\[\bigwedge_n\neg\phi_n.\]
So, we can assume that $|X(T)|$ is uncountable. We will split the proof in various cases given by corresponding lemmas.

\begin{lem}\label{X2NotEmpty} If $X_2(T)\neq\emptyset$, then $T$ is independently axiomatizable.
\begin{proof} In this case there is a sentence, say $\phi_0$, such that $\neg\phi_0$ has continuum many non-isomorphic countable models. By theorem \ref{partitionphi} we know that there are sentences $\{\psialpha|0<\alpha<2^{\aleph_0}\}$ that partition $\neg\phi_0$. Define a new theory \footnote{Note here that both $\psialpha$ and $\phibaralpha$ are defined for $\alpha>0$.} $T'=\{\phibaralpha |0<\alpha<2^{\aleph_0}\}$ by
\[\phibaralpha:\;\; \neg\psialpha\wedge (\neg\phi_0\vee \phialpha).\]

\begin{claim} $T$ and $T'$ are semantically equivalent.
\begin{proof} (of claim)
\[\bigwedge_{\alpha>0} \phibaralpha \leftrightarrow \bigwedge_{\alpha>0} \neg\psialpha\wedge (\neg\phi_0\vee \phialpha) \leftrightarrow
(\bigwedge_{\alpha>0} \neg\psialpha\wedge \neg\phi_0 )\vee
(\bigwedge_{\alpha>0} \neg\psialpha\wedge \phialpha) \leftrightarrow\]
\[(\phi_0\wedge \neg\phi_0 )\vee
(\bigwedge_{\alpha>0} \neg\psialpha\wedge \phialpha) \leftrightarrow(\bigwedge_{\alpha>0} \neg\psialpha\wedge \phialpha) \leftrightarrow
(\bigwedge_{\alpha>0} \neg\psialpha) \wedge (\bigwedge_{\alpha>0}\phialpha) \leftrightarrow (\bigwedge_{\alpha\ge 0}\phialpha).\]
\end{proof}
\end{claim}

\begin{claim} $T'$ is independent.
\begin{proof} (of claim) Let $\alpha<2^{\aleph_0}$ and fix a model $\Malpha$ with $\Malpha\models\psialpha$. By the assumption that the \psialpha's partition $\neg\phi_0$, we get that $\Malpha\models\neg\phi_0$ and for all other $\beta\neq\alpha$, $\Malpha\models\neg\psibeta$. Therefore,
\[\Malpha\models \bigwedge_{\beta\neq\alpha}\phibarbeta\wedge\neg\phibaralpha.\]
This means that $T'\setminus\{\phibaralpha\}\nvDash \phibaralpha$ and $T'$ is independent.
\end{proof}
\end{claim}
Therefore, $T$ is independently axiomatizable.
\end{proof}
\end{lem}

\begin{lem}\label{Step1} If $X_2(T)=\emptyset$ and $|X_0(T)\setminus X_1(T)|=|X(T)|\ge|T_1|$, then $T$ is independently axiomatizable.
\begin{proof}Before we start we need a lemma that essentially is due to Reznikoff (cf. \cite{Reznikoff}) and also appears in \cite{Caicedo}. We include the proof for completeness.
\begin{lem} Let $C,D$ be disjoint sets of sentences with $|D|\le |C|$. If for every $\phi\in C$,
 \[C\cup D\setminus\{\phi\}\nvDash \phi,\]
then $C\cup D$ is independently axiomatizable.
\begin{proof} Let $f:D\rightarrow C$ be a 1-1 function. Then the set \[(C\setminus f(D))\cup\{\phi\wedge f(\phi)|\phi\in D\}\] is independent and semantically equivalent to $C\cup D$.
\end{proof}
\end{lem}

Assume that $|X(T)|=\kappa\ge\omegaone$. By the previous lemma it suffices to find a theory $T'_0$ such that
\begin{itemize}
\item $T'_0\cup T_1$ is equivalent to $T_0\cup T_1$,
\item  $|T'_0|=\kappa\ge|T_1|$ and
\item for every sentence $\phi\in T'_0$, $T'_0\cup T_1\setminus\{\phi\}\nvDash \phi$.
\end{itemize}
Let $T_0=\{\phi_\alpha|\alpha<2^{\aleph_0}\}$ and for every $\alpha$, let $\{\M^{(\alpha)}_n|n\in\mathbb{N}\}$ and $\{\phialphan|n\in\mathbb{N}\}$ be an enumeration of the (countably many) countable models of $\neg\phi_\alpha$ and their Scott sentences respectively. Define
\[\phibaralpha=\bigwedge\{\neg\phialphan|\M^{(\alpha)}_n\notin X_1(T)\mbox{ and }\phialphan\neq\phi^{(\beta)}_m\mbox{ for all }\beta<\alpha,m\in\mathbb{N}\},\]
i.e. we get the conjunction of all the Scott sentences that neither did they appear at a previous step nor their countable model is in $X_1(T)$. If the conjunction is empty we ignore it. By assumption $|X_0(T)\setminus X_1(T)|=|X(T)|=\kappa$ and there have to be $\kappa$ many $\phibaralpha$'s that are not empty.
Let $T'_0=\{\phibaralpha|\alpha<\kappa\}$.

\begin{claim} $T'_0\cup T_1$ is equivalent to $T_0\cup T_1$.
\begin{proof} (of claim) First observe that \[\neg\phialpha\leftrightarrow \bigvee_n\phialphan.\]
Equivalently,
 \[\phialpha\leftrightarrow \bigwedge_n\neg\phialphan.\]
 Thus, $\phialpha\rightarrow\phibaralpha$ and
 \[T_0\cup T_1\models T'_0\cup T_1.\]
 Conversely, let $\M\models T'_0\cup T_1$. We need to prove that $\M\models T_0$, which is equivalent to
 \[\M\models\phialpha, \mbox{ for all $\alpha$,}\]
 which implies that
 \[\M\models\bigwedge_n\neg \phialphan,\mbox{ for all $\alpha$.} \]
 In other words, \[\M\models\neg \phialphan, \mbox{ for all $\alpha$ and $n$}.\]
 Hence, assume that $\M\models \phialphan$, with $\alpha$ minimal with this property and some $n$. Since $\M\models T'_0$, the only case that this can happen is if $\M^{(\alpha)}_n\in X_1(T)$. If $\M$ is not countable, we can pass to a countable elementary submodel (over an appropriate fragment), say $\M_0\prec\M$. Then $\M_0\cong\M^{(\alpha)}_n$ and, therefore, there is $\phi\in T_1$ with $\M_0\models\neg\phi$. If the fragment was chosen to include $\phi$, we also get that $\M\models\neg\phi$, contradicting the fact that $\M\models T_1$.
 \end{proof}
 \end{claim}

 \begin{claim} Every sentence in $T'_0$ is not implied by other sentences in $T'_0\cup T_1$.
\begin{proof} (of claim) Fix $\alpha$ and assume that $\phibaralpha$ is not empty, with, say $\neg\phialphan$, being one sentence in the conjunction. Since $\M^{(\alpha)}_n\models\phialphan$, it cannot satisfy any other Scott sentence and since $\neg\phialphan$ doesn't appear in any other $\phibarbeta$, we conclude that
 \[\M^{(\alpha)}_n\models\neg\phibaralpha\wedge\bigwedge_{\beta\neq\alpha}\phibarbeta.\]
 So, $\M^{(\alpha)}_n\models T'_0\setminus\{\phibaralpha\}$. But also, for every $\phi\in T_1$, $\M^{(\alpha)}_n\models\phi$, because otherwise we would have $\M^{(\alpha)}_n\in X_1$ and this would prevent $\neg\phialphan$ from being in the conjunction of $\phibaralpha$. Contradiction. Putting everything together we get that
 \[\M^{(\alpha)}_n\models T'_0\cup T_1\setminus\{\phibaralpha\},\]
 witnessing that $T'_0\cup T_1\setminus\{\phibaralpha\}\nvDash \phibaralpha$.
 \end{proof}
 \end{claim}

 This finishes the proof.
 \end{proof}
 \end{lem}

 \begin{lem} If $X_2=\emptyset$ and $|X(T)|\ge|T_1|$, then $T$ is independently axiomatizable.
 \begin{proof} If $|X(T)|=\kappa$ and $|X_1(T)|<\kappa$, then the assumptions of lemma \ref{Step1} are satisfied and $T$ is independently axiomatizable. So, assume that $X_1(T)=\{\Malpha|\alpha<\kappa\}$ with $\kappa\ge\omegaone$ and $T_1=\{\psi_\alpha|\alpha<\lambda\}$ with $\lambda\le\kappa$. We will find another theory $\Tstar$, equivalent to $T$ and for which
 \[|X_0(\Tstar)\setminus X_1(\Tstar)|=|X(\Tstar)|=|X(T)|\ge|T_1|=|\Tstar_1|.\]
  Again, by lemma \ref{Step1} we are done. We know that the only case that a sentence $\phi$ can be in $T_1$ is if for all countable $\alpha$, both $\Phi_\alpha(\neg\phi)$ and $\Psi_\alpha(\neg\phi)$ are countable. For every $\alpha<\omegaone$ define new sets $C_\alpha(\neg\phi)$ and $S_\alpha(\neg\phi)$:

  $\phialphapar{\emptyset,\M}\in C_\alpha(\neg\phi)$ if and only if $\phialphapar{\emptyset,\M}\in\Phi_\alpha(\neg\phi)$ and there are only countably many countable models of $\phi$ that satisfy $\phialphapar{\emptyset,\M}$, and

  $\sigma\in S_\alpha(\neg\phi)$ if and only if there exists a countable model $\M$ that satisfies some $\phialphapar{\emptyset,\M}\in C_\alpha(\neg\phi)$ and $\sigma$ is its Scott sentence.

  Both $C_\alpha(\neg\phi)$ and $S_\alpha(\neg\phi)$ are countable for all $\alpha$, since $\Phi_\alpha(\neg\phi)$ is countable. We will distinguish three cases:

     \textbf{Case I:} $\kappa>\omegaone$ and $cf(\kappa)\neq\omegaone$. Since\[\kappa=|X_1(T)|\le|T_1|\aleph_1=\lambda\cdot\aleph_1\] and \[\kappa>\omegaone,\] it must be that $\kappa=\lambda$. Since $cf(\kappa)\neq\omegaone$, there exists an ordinal $\gamma<\omegaone$ and $\kappa$ non-isomorphic countable models in $X_1$ of Scott height less than $\gamma$. Define inductively a new theory, considering the sentence $\neg\psi_\alpha$ at stage $\alpha<\lambda$. Choose $\beta$ larger than $\gamma$ and replace $\psi_\alpha$ by \[\psi_\alpha^{(0)}:=\bigwedge \{\neg\sigma|\sigma\in S_\beta(\neg\psi_\alpha) \}\] and
     \[\psi_\alpha^{(1)}:=\psi_\alpha\vee\bigvee\{\sigma|\sigma\in S_\beta(\neg\psi_\alpha) \}.\]
     It is not hard to see that $\psi_\alpha$ is equivalent to the conjunction of $\psi_\alpha^{(0)}$ and $\psi_\alpha^{(1)}$. Also, observe that $\neg\psi_\alpha^{(0)}$ has countably many countable models and all the countable models of $\neg\psi_\alpha$ of Scott height less than $\gamma$ satisfy it. Repeating this for $\lambda$ many steps we will get eventually a theory $\Tstar$ such that $X_0(\Tstar)$ will contain all countable models that are in $X_1(T)$ that have Scott height $<\gamma$. By the assumption on $\gamma$, \[|X_0(\Tstar)\setminus X_1(\Tstar)|=\kappa.\]

     \textbf{Case II:} $\kappa>\omegaone$ and $cf(\kappa)=\omegaone$. As before $\kappa=\lambda$, but the difference now is that we may not have an ordinal $\gamma$ as before. Instead, assume that there are cardinals $\{\mu_i|i<\omegaone\}$ and countable ordinals $\{\alpha_i|i<\omegaone\}$ such that
     \begin{itemize}
     \item for all $i<j$, $\omegaone<\mu_i<\mu_j$,
     \item $\sup_i \mu_i=\kappa$,
     \item for $i=0$, $\alpha_0=0$,
     \item for all $i<j$, $\alpha_i<\alpha_j$, and
     \item for $j$ limit    ordinal, $\sup_{i<j} \alpha_i=\alpha_j$,    and
     \item for all $i<\omegaone$, the number of countable models in $X_1(T)$ that have Scott    height $\alpha$ with $\alpha_i\le\alpha<\alpha_{i+1}$ is equal to    $\mu_i$.
     \end{itemize}
     This also yields a partition $T_1=\bigcup_{i<\omegaone} T_1^{(i)}$ such that for all $i$
     \begin{itemize}
     \item for all $\psi\in T_1^{(i)}$, $\neg\psi$ has a countable model    of Scott height $\alpha$, $\alpha_i\le\alpha<\alpha_{i+1}$  and
         \item $|\{\M|\M\models\neg\psi,\mbox{ some }\psi\in T_1^{(i)},\mbox{ $\M$ countable and }    \alpha_i\le\scM<\alpha_{i+1}\}|    =\mu_i$.
     \end{itemize}
     As before we define a new theory inductively: At stage $\alpha$, if $\psi_\alpha\in T_1^{(i)}$, choose $\beta\ge \alpha_{i+1}$ and replace $\psi_\alpha$ by
     \[\psi_\alpha^{(0)}:=\bigwedge \{\neg\sigma|\sigma\in S_\beta(\neg\psi_\alpha) \}\]
     and
     \[\psi_\alpha^{(1)}:=\psi_\alpha\vee\bigvee\{\sigma|\sigma\in S_\beta(\neg\psi_\alpha) \}.\]
     It is not hard to see that $\psi_\alpha$ is equivalent to the conjunction of $\psi_\alpha^{(0)}$ and $\psi_\alpha^{(1)}$. Also, $\neg\psi_\alpha^{(0)}$ has countably many countable models, while $\neg\psi_\alpha^{(1)}$ has $\alephs{1}$ many countable models, and all the countable models of $\neg\psi_\alpha$ of Scott height $<\alpha_{i+1}\le\beta$ satisfy $\neg\psi_\alpha^{(0)}$. Eventually, after $\lambda$ many steps we will get a theory $\Tstar$ such that $X_0(\Tstar)\setminus X_1(\Tstar)$ contains at least $\mu_i$ many countable models $\M\in X_1(T)$ that have Scott height $\alpha_i\le\scM<\alpha_{i+1}$. By the assumptions on the $\mu_i$'s,
     \[|X_0(\Tstar)\setminus X_1(\Tstar)|=\kappa.\]

     \textbf{Case III:} $\kappa\le\omegaone$. Then $\lambda\le\omegaone$ and we can use Caicedo's theorem (in \cite{Caicedo}) that every set with $\le\omegaone$ sentences in $\lomegaone$ is independently axiomatizable.
\end{proof}
\end{lem}

\begin{thrm}\label{MainTheorem} If $\lang{}$ is a countable language, $T$ a theory in $\lomegaone$, $T_1$ is as in definition \ref{t0t1t2def} and $T_1$ is also small in $T$ (cf. definition \ref{SmallDef}), then $T$ is independently axiomatizable.
\begin{proof} If $X_2\neq\emptyset$, then use lemma \ref{X2NotEmpty}. If it is empty, then use the previous lemma.
\end{proof}
\end{thrm}

\begin{cor}\label{VCFailsCor} If the Vaught Conjecture holds, then every $T\subset\lomegaone$ is independently axiomatizable.
\begin{proof} The Vaught Conjecture gives us that $T_1=\emptyset$. Then use the previous theorem.
\end{proof}
\end{cor}

\begin{cor}If $|X(T)|=\continuum$, then $T$ is independently axiomatizable.
\begin{proof} Then $|X(T)|\ge|T_1|$ and we can again apply theorem \ref{MainTheorem}.
\end{proof}
\end{cor}

\section{Reformulations and open questions}

\hspace{.5cm} In this section we reformulate the previous theorems as statements about Borel sets and give some open problems. Recall that a collection of Borel sets $\B=\{B_i|i\in I\}$ is independent if $\bigcap\B\neq\emptyset$ and for every $i\in I$, $\bigcap_{j\neq i} B_j\setminus B_i\neq\emptyset$, and that two collections $\B,\B'$ are equivalent if $\bigcap \B=\bigcap\B'$.

\begin{thrm}\label{MainTheorem2} Every collection of Borel sets $\B=\{B_i|i\in\continuum\}$ with $\bigcap\B\neq\emptyset$ admits an equivalent independent collection.
\begin{proof}The proof closely resembles the proofs of lemma \ref{Step1} and lemma \ref{X2NotEmpty}. We have two cases:

\textbf{Case I:} There is an $i_0\in I$, such that $\complement B_{i_0}$, the complement of $B_{i_0}$, is uncountable. Then we can partition $\complement B_{i_0}$ into continuum many sets $\bigcup_{x\in \complement B_{i_0}} \{x\}$. Call these sets $\{C_j|j\neq i_0,j<\continuum\}$. Define now a new collection of Borel sets $\B'=\{B'_j|j\neq i_0,j<\continuum\}$ by
\[\B'_j:= \complement C_j\cap (\complement B_{i_0}\cup B_j).\]

\begin{claim} $\B$ and $\B'$ are equivalent.
\begin{proof} (of claim)
\[\bigcap_{j\neq i_0} B'_j=
\bigcap_{j\neq i_0} \complement C_j\cap (\complement B_{i_0}\cup B_j)=
(\bigcap_{j\neq i_0} \complement C_j\cap \complement B_{i_0} )\cup
(\bigcap_{j\neq i_0} \complement C_j\cap B_j)=\]
\[(B_{i_0} \cap \complement B_{i_0} )\cup
(\bigcap_{j\neq i_0} \complement C_j\cap B_j)=
\bigcap_{j\neq i_0} \complement C_j\cap B_j=
\bigcap_{j\neq i_0} \complement C_j\cap \bigcap_{j\neq i_0} B_j=
\bigcap_{j} B_j.\]
\end{proof}
\end{claim}

\begin{claim} $\B'$ is independent.
\begin{proof} (of claim) Let $x\in C_j$. By the properties of the $C_j$'s, we get that $x\in \complement B_{i_0}$ and $x\notin C_{j'}$, for $j'\neq j$. Therefore,  $x\notin B'_j$, while $x\in B'_{j'}$, $j'\neq j$, which implies that
\[\bigcap_{j'\neq j} B'_{j'}\setminus B'_j\neq\emptyset.\]
\end{proof}
\end{claim}

\textbf{Case II:} For all $i\in I$, $\complement B_{i}$, the complement of $B_{i}$, is countable. Construct a new collection $\B'=\{B'_j|j<2^{\aleph_0}\}$ with
\[B'_j=B_j \bigcup_{i<j}\complement B_i.\]
If the set is equal to the whole space, we ignore it and proceed to the next one. Observe here that the the complement of $B'_j$ is a subset of the complement of $B_j$, which is countable by assumption. Therefore, it is Borel.
\begin{claim} $\B$ and $\B'$ are equivalent.
\begin{proof} (of claim) It is immediate that $B'_j\supset B_j$, which means that $\bigcap_j B'_j\supset \bigcap_j B_j$. So, let $x\in \bigcap_jB'_j$. By induction on $j$ we can prove that $x\in B_j$. Assume that $x\in \bigcap_{i<j} B_i$. Then $x\notin\bigcup_{i<j}\complement B_i$. Since, $x\in B'_j$, this implies that $x\in B_j$. Consequently, $x\in \bigcap_j B_j$.
\end{proof}
\end{claim}

\begin{claim} $\B'$ is independent.
\begin{proof} (of claim) Fix $j<2^{\aleph_0}$ and assume that $B'_j$ is not equal to the whole space. Say $y\in \complement B'_j$ witnesses this. Then,  $y\in\complement B_j$ and by definition $y\in B'_i$, for all $i>j$. Similarly, $y\notin B'_j$ implies that $y\notin\bigcup_{i<j}\complement B_i$, which means that $y\in\bigcap_{i<j} B_i$. Then, $y\in B'_i$, for all $i<j$, and over all, $y\in \bigcap_{i\neq j} B'_i\setminus B'_j$.
\end{proof}
\end{claim}

In either case, we constructed an independent collection of Borel sets $\B'$ which is equivalent to $\B$.
\end{proof}
\end{thrm}

It would be interesting if we could derive theorem \ref{MainTheorem} from theorem \ref{MainTheorem2}. This would eliminate the extra assumptions of theorem \ref{MainTheorem}.

\begin{df} Let $T\models_g\phi$ mean that in all generic extensions every model of $T$ is also a model of $\phi$.
\end{df}

This is a stronger notion than $T\models\phi$ and is related to $T\vdash_{\lomegaone}\phi$, but we will not define $\vdash_{\lomegaone}$ here. We can prove
\begin{thrm}\label{Models_g}If $T\models_g\phi$, then there are countably many sentences $\phi_0,\phi_1,\ldots\in T$ such that
\[\bigwedge_n \phi_n\models_g\phi.\]
\end{thrm}

We now ask whether we can replace $\models$ by $\models_g$ in theorem \ref{MainTheorem}. The problem is that $T$ and $T'$ may not be semantically equivalent in a generic extension. This is an open question we did not consider. We can also reformulate this problem using the language of Boolean Algebras. We know that the $\lomegaone$- sentences form a $\sigma$-complete Boolean Algebra with $\phi\le\psi$  if and only if $\phi\rightarrow\psi$. Using theorem \ref{Models_g} we can prove that the $\sigma$-filter generated by a theory $T$ is equal to
\[T'=\{\psi|T\models_g\psi\}.\]

\begin{df} A set $A$ of sentences is called $\sigma$-\emph{filter independent}, if for all $\phi$, $\phi$ is not in the $\sigma$-filter generated by $A\setminus\{\phi\}$.
\end{df}
The problem is given a set of sentences $A$ to find another set$A'$ such that
\begin{itemize}
\item $A$ and $A'$ generate the same $\sigma$-filter and
\item $A'$ is $\sigma$-filter independent.
\end{itemize}

We can also extend the question to finding conditions under which a Boolean Algebra satisfies the above statement. As far as we know this problem is open. Another extension would be to prove that any $\lomegaone$ theory is independently axiomatizable, without assuming countability of the language. Our techniques here rely heavily on this assumption.

\section{changes}
\begin{itemize}
  \item Definition 10: The line that starts "For a $\lomegaone$..., a space was added between "and" and $\alpha$.
    \item Definition 10: The line that starts "This enables us to consider..., a space was added between "and" and $\Phi_{\gamma+1}$.
  \item Definition 11: The set \[\{M\in Mod(\lang{}): M\models \varphi(n_1,...n_m)\}\] became \[\{M\in Mod(\lang{})| M\models \varphi(n_1,...n_m)\},\] i.e. the symbol ``$|$" replaced the symbol ``$:$"
  \item Proof of Lemma 15: The set   \[\{(x,\M)\in P\times Mod(\phi)|\M\models t(x),\mbox{ Scott height}(\M)<\beta\}.\] became \[\{(x,\M)\in P\times Mod(\phi)|\M\models t(x),\alpha(\M)<\beta\},\] i.e. $\alpha(\M)$ replaced ``Scott height$(\M)$". 
  \item Definition 21: The conjuction \[\bigwedge\neg\phi_n.\] became \[\bigwedge_n\neg\phi_n,\] i.e. a subscript $n$ was added.
  \item Proof of Lemma 25, Case III: A space was added in the citation of Caicedo's theorem. So, ``... we can use Caicedo's theorem (in\cite{Caicedo})" became ``we can use Caicedo's theorem (in \cite{Caicedo})".
  \item The statement of Theorem 26 was rewritten.
  \item Theorem 29, Claim 9: In the line that starts "Similarly, $y\notin B'_j...$, the next three instances of the variable $x$ were replaced with the variable $y$. Also, the disjunction \[y\in \bigcup_{i\neq j} B'_i\setminus B'_j\] was replaced with a conjunction \[y\in \bigcap_{i\neq j} B'_i\setminus B'_j.\]
  \item Before Definition 30, in the sentence that starts ``In either case...", the word ``which" replaced the word ``with".
\end{itemize}


\begin{thebibliography}{99}
\bibitem{Kechrisbook} Alexander S. Kechris, \emph{Classical Descriptive SetTheory}, Graduate Texts in Mathematics, 156, Springer-Verlag, NewYork, 1995.
\bibitem{Caicedo} X. Caicedo, \emph{Independent Sets of Axioms in$L_{\kappa\alpha}$}, Canad.Math.Bull., Vol.\textbf{24} (2), 1981,pp.219-223
\bibitem{Reznikoff} M. I. Reznikoff, \emph{Tout ensemble de formulesde la logique classique est equivalent $\acute{a}$ un ensembleindependant}, C.R. Acad. Sc. Paris, \textbf{260}, 2385-2388(1965).
\bibitem{BeckerKechrisBook} H. Becker \& A. Kechirs, \emph{TheDescriptive Set Theory of Polish Group Actions}, CambridgeUniversity Press, London Mathematical Society Lecture Note Series232, 1996.
\bibitem{Miller'sWebpage} http://www.math.wisc.edu/~miller/res/problem.pdfThis webpage contains a list of intersting problems in Set Theoryand Model Theory.

\end{thebibliography}
\end{document}